\documentclass[11pt,a4paper]{article}
\usepackage{epsfig}
\usepackage{amssymb}
\usepackage{amsmath}
\usepackage{epic}
\usepackage{eepic}
\usepackage{enumerate}

\def\T{{\mathbb T}}

\def\cT{{\cal T}}
\def\U{{\cal U}}
\def\V{{\cal V}}

\def\P{{\mathbb P}}
\newtheorem{theorem}{Theorem}[section]
\newtheorem{lemma}[theorem]{Lemma}

\newtheorem{example}[theorem]{Example}
\newtheorem{problem}[theorem]{Open Problem}

\author{Erik I. Broman\footnote{Department of Mathematics, Uppsala University,
url: \texttt{http://www2.math.uu.se/\~{ }broman/}, e-mail:broman@math.uu.se
supported by the Swedish Research Council.}}
\date{\today}

\title{Stochastic Ordering of Infinite Binomial Galton-Watson Trees.}

\begin{document}
  \maketitle

  \begin{abstract}
We consider Galton-Watson trees with ${\rm Bin}(d,p)$ offspring distribution. We
let $T_{\infty}(p)$ denote such a tree conditioned on being infinite.
For $d=2,3$ and any $1/d\leq p_1 <p_2 \leq 1$, we show that there exists a coupling
between $T_{\infty}(p_1)$ and $T_{\infty}(p_2)$ such that $\P(T_{\infty}(p_1) \subseteq T_{\infty}(p_2))=1.$




  \end{abstract}

\section{Introduction} \label{secintro}
We start with a somewhat informal motivation of the paper, and give formal definitions below. For any $d\geq 2$
and $0<p<1$
we let $T(p)=T(p,d)$ denote a Galton-Watson tree, with binomial offspring distribution, with parameters $p,d$. For fixed $d\geq 2$ and
$0<p_1<p_2<1,$
a trivial coupling argument allows us to couple the trees $T(p_1)$ and $T(p_2)$ such that $\P(T(p_1) \subseteq T(p_2))=1$.
A natural question to ask is whether this property is preserved if we condition the trees on being infinite; if we let
 $T_\infty(p)$
denote a sample of $T(p)$ conditioned on being infinite, is it the case that there exists a coupling of
$T_\infty(p_1)$ and  $T_\infty(p_2),$ where $1/d\leq p_1<p_2<1,$ such that
\begin{equation}\label{eqn:conj}
\P(T_\infty(p_1) \subseteq T_\infty(p_2))=1? 
\end{equation}
The case $p_1=1/d$ warrants an explanation given below.

A finite version of this was proved by Luczak and Winkler in \cite{LW} (see Theorem \ref{thmfinite} below for a 
precise statement of their result). Their result will be a key ingredient of the proof of our main result.
The analogue of \eqref{eqn:conj} for Poisson offspring distribution was proven by Lyons, Peled and Schramm
in \cite{LPS}. It is natural to ask whether such a result would hold for any (parametrized) offspring distribution
for which the unconditioned trees can be appropriately coupled. Example \ref{ex:counter} shows that the answer is no.
This counterexample is a variant of the one by Janson in \cite{Janson}, used to prove that the finite
version of \eqref{eqn:conj} does not hold for general offspring distributions
(see also the remark after Theorem \ref{thmfinite}).

\begin{example} \label{ex:counter}
{\rm Consider a Galton Watson tree $T(r)$ with the following offspring distribution. Any individual in the tree has
$0,1$ or $2$ children with probability $1/2-2r,r$ and $1/2+r$ respectively, where $r\in[0,1/4].$
It is not hard to see that for $0\leq r_1<r_2\leq 1/4,$ we can couple the construction of $T(r_1)$
and $T(r_2)$ such that $\P(T(r_1)\subseteq T(r_2))=1.$ Furthermore, if we let $A$ be the event that the root
has exactly one offspring, we see that 
\[
\P(A | |T(r)|=\infty)=\frac{\P(|T(r)|=\infty | A) \P(A)}{\P(|T(r)|=\infty)}=\P(A)=r.
\]
Thus, for $r_1<r_2$, the probability that the root has only one offspring is larger for $T(r_2)$ 
conditioned on survival than for $T(r_1)$ conditioned on survival. 
Therefore, \eqref{eqn:conj} cannot hold for this offspring distribution.}
\end{example}

Similar questions can be asked about percolation clusters on graphs. One of the most interesting cases is
 the corresponding problem for bond-percolation on ${\mathbb Z}^d$; see Open problem
\ref{q3} at the very end of the paper. In this paper we will study the case of Galton-Watson trees with binomial offspring
for $d=2,3$. It is an open problem whether this result holds also for $d \geq 4$ (see Section \ref{secopenprob}).

We proceed to give some definitions needed for the statements of the main results.
Let $\T^d$ be the rooted, ordered, labeled tree, in which every vertex including the root has exactly $d$ children
ordered from left to right.
Hence, $\T^d$ is a regular tree in which every vertex has degree $d+1$ except for the root,
which has degree $d$. We will let $o$ denote the root, and the labeling of $\T^d$ is done in the natural
way, so that we identify the vertex set $V(\T^d)$ with the set $\{o\}\bigcup_{n=1}^{\infty} \{1,\ldots,d\}^n$.
For any two elements $u=(u_1,\ldots,u_k),v=(v_1,\ldots,v_l)\in V(\T^d)\setminus\{o\}$
we let $(u,v)=(u_1,\ldots,u_k,v_1,\ldots,v_l)$
denote the concatenation of $u$ and $v.$ For $i\in \{1,\ldots,d\},$ we will allow a slight abuse of notation, by writing
$(u,i)$ instead of $(u,(i))$ for $(u_1,\ldots,u_k,i).$
Furthermore, for any $u\in V(\T^d)$
we let $(u,o)=(o,u)=u$ so that in particular, $(o,o)=o$. 
We will use the natural edge set for $\T^d$ and all other trees, but since the edge set is determined
by the vertex set, it will not play any part in the analysis. We will therefore only refer to a tree by
its vertex set.

A {\em subtree} of $\T^d$ is defined to be a connected subgraph
of $\T^d.$ For any such subtree $T,$ we will let $V(T)$ denote the vertex set. Furthermore, we let $|T|$ denote the
number of vertices of $T,$ and we call this the {\em size} of $T$. We let
$c_k=c_k(d)$ denote the number of subtrees $T$ such that $o\in V(T)$ and $|T|=k.$
For example, $c_3(2)=5$ and $c_2(d)=d.$
For $u\in V(\T^d)$ and a subtree $T,$ let $T^u$ denote the subtree
(of $T$) with vertex set $V(T^u):=\{v\in V(T): v=(u,w) \textrm{ for some } w\in V(\T^d)\}.$ Note that if $u\not \in V(T),$
we get that $T^u=\emptyset.$ Informally, $T^u$ is simply the tree consisting of $u$ and the descendants of $u$
that belongs to $T.$
 We also define $H(T^u):=\{w\in V(\T^d): (u,w)\in V(T^u)\},$ which is simply a shift of $T^u,$
mapping $u$ to $o.$ For $i\in \{1,\ldots,d\}$ we sometimes abuse notation and write $T^i$ instead of $T^{(i)}.$

For $d\geq 2$ and $0<p<1,$ perform site percolation (see \cite{Grimmett} for a general overview on percolation)
with density $p$ on $\T^d$. We consider the resulting random subgraph of $\T^d$, and let $T(p)$ be the component of the root.
If the root is removed in the percolation procedure, we take $T(p)=\emptyset.$
It is clear that $T(p)$ is the family tree of a Galton-Watson process with a Bin$(d,p)$
offspring distribution (see \cite{AN} for a general overview on Galton-Watson processes), except that here, we allow
this family tree to be the empty set, which is only a matter of convenience. We denote the distribution of $T(p)$
by $\cT(p).$
Define $\eta_k(p):={\mathbb P}(|T(p)|=k)$ for $0\leq k \leq \infty$,
(note that we include the case $k=\infty$). It is not hard to check that
\begin{equation} \label{eqn10}
\eta_k(p)=c_k p^k(1-p)^{(d-1)k+1},
\end{equation}
if $k<\infty$. For $1\leq k<\infty,$ let $T_k$ be uniformly chosen among the subtrees $T,$ such that $o\in V(T)$ and $|T|=k.$
Let the distribution of $T_k$ be denoted by $\cT_k.$ It is not hard to check that
the distribution of $T(p),$ conditioned on the event that $|T(p)|=k,$ is also $\cT_k$ (in particular, it is
independent of $p$). Sometimes, it will be convenient to think of the empty set as the tree of size $0,$
and then we will use the notation $T_0=\emptyset.$

It is well known, that for $p > 1/d$, $\P(|T(p)|=\infty)>0$. For such $p>1/d,$
let $T_{\infty}(p)$ denote a random subtree of $\T^d$ whose distribution equals that of $T(p),$ conditioned on the event
$|T(p)|=\infty.$ It is also well known that $\P(|T(1/d)|=\infty)=0.$ However, one can still define an infinite critical random tree
$T_{\infty}(1/d)$ in a natural way. This random tree is the so-called incipient infinite cluster 
(see \cite{Kesten} or \cite{Angel}) on $\T^d.$  In order to define this, 
let $(Z_i)_{i\geq 1}$ be an i.i.d.\ sequence such that $\P(Z_i=j)=1/d$ for $j\in\{1,\ldots,d\}.$ Furthermore, let
$(T_{i,j}(1/d))_{i\geq 1, 1\leq j \leq d}$ be an i.i.d.\ sequence, independent of $(Z_i)_{i\geq 1},$ and such that
$T_{i,j}(1/d) \sim \cT(1/d)$ for every $i,j.$ Informally, we will let $T_{\infty}(1/d)$ be the tree consisting
of a single infinite line (backbone) determined by the sequence $(Z_i)_{i\geq 1},$ and onto this line we attach the trees
$(T_{i,j}(1/d))_{i\geq 1, 1\leq j \leq d}$ in the appropriate places.
Formally, we let
\[
V(T_{\infty}(1/d)):=\{o\}\bigcup_{i=1}^\infty \{(Z_1,\ldots,Z_i)\} \bigcup_{j\in \{1,\ldots,d\}\setminus \{Z_i\}}\bigcup_{u\in T_{i,j}(1/d)}\{(Z_1,\ldots,Z_{i-1},j,u)\},
\]
where for $i=1,$ we let $(Z_1,\ldots,Z_{i-1},j,u)=(j,u).$
Observe that if for some $i,j$ we have $T_{i,j}(1/d)=\emptyset,$ then this will not make any contribution to the vertex set
of $T_{\infty}(1/d)$.
For $p\geq 1/d,$ we denote the distribution of $T_\infty(p)$ by $\cT_\infty(p).$

For any two subtrees $S, T$ of $\T^d$, we write $S \subseteq T$ if $S$ is a subgraph of $T$.
As mentioned above, in \cite{LW}, the following theorem was proved, which we restate here since it will
be crucial for proving our main results.
\begin{theorem}[Luczak, Winkler] \label{thmfinite}
For any $d\geq 2,$ there exists a coupling of $(T_k)_{k\geq 0}$
(where $T_k \sim \cT_k$ for every $0\leq k<\infty$)
such that
\begin{equation} \label{eqn1}
\P(T_0 \subset T_1 \subset T_2 \subset \cdots \subset T_{l} \subset \cdots)=1.
\end{equation}
Furthermore, using this coupling, we have that
\begin{equation} \label{eqn3}
\bigcup_{k=0}^\infty T_k \sim \cT_{_\infty}(1/d).
\end{equation}
\end{theorem}
{\bf Remark:} 
As mentioned above, it is proved in \cite{Janson} that equation (\ref{eqn1}) does not hold for 
general offspring distributions.

Theorem \ref{thmmain} is natural in light of Theorem \ref{thmfinite} and is the main result of this paper.
\begin{theorem} \label{thmmain}
For $d=2,3$ and $1/d\leq p_1<p_2\leq 1$ there exists a coupling of $T_{\infty}(p_1)$ and $T_{\infty}(p_2)$
(where $T_{\infty}(p_1)\sim \cT_{\infty}(p_1)$
and $T_{\infty}(p_2)\sim \cT_{\infty}(p_2)$) such that
\[
\P(T_{\infty}(p_1) \subset T_{\infty}(p_2))=1.
\]
\end{theorem}
{\bf Remark:} As mentioned before, the corresponding result for  Galton-Watson trees with
Poisson offspring distributions was proved in \cite{LPS}.

\medskip

We also prove the following very natural theorem.
\begin{theorem} \label{thmmain2}
For $d=2,3$, $1/d\leq p \leq 1$ and any $k$, there exists a coupling of $T_k$ and $T_{\infty}(p),$ (where
$T_k\sim \cT_k$ and $T_{\infty}(p)\sim \cT_{\infty}(p)$ )
such that
\[
\P(T_k \subset T_{\infty}(p))=1.
\]
\end{theorem}
{\bf Remark:} Of course Theorem \ref{thmmain2} is a trivial corollary of Theorems \ref{thmfinite}
and \ref{thmmain}. However, in $d=2$, we first prove Theorem \ref{thmmain2} and use it to prove
Theorem \ref{thmmain}. When $d=3,$ the main effort will be to prove Theorem \ref{thmmain} for the special 
case $p_1=1/3,$ from which Theorem \ref{thmmain2} then follows. In turn, Theorem \ref{thmmain2} will then
be used to prove Theorem \ref{thmmain} for every $1/3<p_1<p_2\leq 1.$

\medskip
The outline of the rest of the paper is as follows.
All results for $d=2$ are proved in Section \ref{secd=2} while the results for
$d=3$ are proved in Section \ref{secd=3}. In Section \ref{secopenprob}, we present some open problems.

\section{The case $d=2$} \label{secd=2}

We start with a preliminary result which will be useful to us. In \cite{LW}, it is proved that
\[
c_k(d)=\frac{{dk \choose k}}{(d-1)k+1}.
\]
The following lemma is an easy consequence of this, and is therefore left without proof.

\begin{lemma} \label{lemmagen1}
For any $d$, the sequence $(c_{k-1}/c_k)_{k\geq 1}$ decreases in $k$, and furthermore
\[
\lim_{k \rightarrow \infty}\frac{c_{k-1}}{c_k}=\frac{1}{d}\left(\frac{d-1}{d}\right)^{d-1}.
\]
\end{lemma}
{\bf Remark:} By maximizing $p(1-p)^{d-1}$ (with respect to $p$), it follows from \eqref{eqn10} that
\[
\frac{\P(|T(p)|=k)}{\P(|T(p)|=k-1)}=\frac{c_{k}}{c_{k-1}}p(1-p)^{d-1} \leq 1,
\]
so that $\P(|T(p)|=k)$ is decreasing in $k$. We conclude, for future reference,
that for any $l \leq k$ and any $p$,
\begin{equation} \label{eqn5}
\frac{c_{k-l}}{c_k}\geq \left(\frac{1}{d}\left(\frac{d-1}{d}\right)^{d-1}\right)^l
\geq \left(p(1-p)^{d-1}\right)^l.
\end{equation}

\medskip

We assume throughout the rest of this section that $d=2$ and $p \geq 1/2.$ 
Furthermore, any tree in the rest of this section will be a subtree of $\T^2$.

Informally the main idea is as follows. Consider a tree $ T_\infty \sim \cT_\infty(p)$ and the 
two subtrees $T^1_\infty, T^2_\infty.$ One of these will necessarily be infinite, while the 
other may be finite. Thus, one way of generating $ T_\infty$ should be to start with a root o, 
then to pick one of the children $(1),(2)$ with equal probability and attach
 an independent copy from 
$\cT_\infty(p)$ to it. Then, we use another random tree $T^*(p)$ (with a very particular distribution),
and attach this tree to the second child. This is made precise in Lemma \ref{lemmad=21}. 
This lemma can then be used to prove Theorem \ref{thmmain2}. 
In Lemma \ref{lemmad=22}, we prove that for any $1/2\leq p_1<p_2\leq 1,$ 
we can take $|T^*(p_1)|\leq |T^*(p_2)|$, which 
will allow us to prove Theorem \ref{thmmain}.

We start by constructing $T^*(p).$
To that end, let $(T_k)_{k\geq 1}$ be a sequence of random trees such that $T_k \sim \cT_k$ for every $k,$
and let $T_\infty(p) \sim \cT_\infty(p).$ Furthermore, let $U\sim U([0,1])$ be independent of
$(T_k)_{k\geq 1}$ and $T_\infty(p).$ If $U\leq 2p\eta_{0}(p),$ we let $T^*(p)=\emptyset,$ while if
$\sum_{l=0}^{k-1}2p\eta_{l}(p)<U\leq \sum_{l=0}^{k}2p\eta_{l}(p)$
for some $1\leq k<\infty,$ we let $T^*(p)=T_k,$ and otherwise we let $T^*(p)=T_\infty(p).$ We observe that
$\P(|T^*(p)|=k)=2p\eta_k(p)$ for $k<\infty,$ and that $\P(|T^*(p)|=\infty)=p \eta_\infty(p).$ Indeed,
an elementary argument shows that
$\eta_\infty(p)=p(1-(1-\eta_\infty(p))^2),$ so that $\eta_\infty(p)=(2p-1)/p$ when $p>1/2,$ while $\eta_\infty(1/2)=0.$
Therefore we have that
$p\eta_\infty(p)+\sum_{k=0}^\infty2p\eta_k(p)=p\eta_\infty(p)+2p(1-\eta_\infty(p))=1.$
Note that if
$|T^*(p)|>0,$ we have that $o\in V(T^*(p)).$ We denote the distribution of
$T^*(p)$ by $\cT^*(p).$

We can now prove the following easy lemma.
\begin{lemma} \label{lemmad=21}
For $p\geq 1/2$, let $(X,T_\infty(p),T^*(p))$ be three independent random variables, 
where $T_\infty(p)\sim \cT_\infty(p),T^*(p)\sim \cT^*(p)$
and $X\in\{1,2\}$ is such that ${\mathbb P}(X=1)={\mathbb P}(X=2)=1/2.$
Define the tree $\tilde{T}(p)$ by letting
\begin{equation} \label{eqn6}
V(\tilde{T}(p))=\{o\}\bigcup_{u\in V(T_\infty(p))}\{(X,u)\} \bigcup_{v\in V(T^*(p))}\{(3-X,v)\}.
\end{equation}
We have that $\tilde{T}(p)\sim \cT_{\infty}(p).$
\end{lemma}
{\bf Remark.}
Thus, $\tilde{T}(p)$ is constructed by starting with a root, and then attaching the tree $T_\infty(p)$ either to the
left or to the right of the root depending on the value of $X,$ and then attaching $T^*(p)$ to the other side.
\medskip

\noindent {\bf Proof.}
Recalling the notation of Section \ref{secintro}, we see that for $k<\infty,$
\begin{eqnarray} \label{eqn4}
\lefteqn{{\mathbb P}(|T^1(p)|=k,|T^2(p)|=\infty \, |\, |T(p)|=\infty)=\P(|T^1(p)|=\infty,|T^2(p)|=k  \, |\, |T(p)|=\infty)} \\
& & =\frac{p \eta_{\infty}(p) \eta_k(p)}{\eta_{\infty}(p)}=p\eta_k(p)
=\frac{1}{2}2p\eta_k(p)=\P(X=2, |T^*(p)|=k, |T_\infty(p)|=\infty) \nonumber \\
& & ={\mathbb P}(|\tilde{T}^1(p)|=\infty,|\tilde{T}^2(p)|=k)
={\mathbb P}(|\tilde{T}^1(p)|=k,|\tilde{T}^2(p)|=\infty), \nonumber
\end{eqnarray}
since obviously $\P(|T_\infty(p)|=\infty)=1.$
We also have that for $p>1/2,$
\begin{eqnarray*} 
\lefteqn{{\mathbb P}( |T^1(p)|=\infty, |T^2(p)|=\infty \, |\, |T(p)|=\infty)=\frac{p \eta_{\infty}(p) \eta_{\infty}(p)}{\eta_{\infty}(p)}=p\eta_{\infty}(p)} \\
& & =\P(|T^*(p)|=\infty, |T_\infty(p)|=\infty)
={\mathbb P}(|\tilde{T}^1(p)|=\infty,|\tilde{T}^2(p)|=\infty). \nonumber
\end{eqnarray*}
Note that in the case $p=1/2,$ only (\ref{eqn4}) is relevant.
We conclude that $(|\tilde{T}^1(p)|,|\tilde{T}^2(p)|)$ and $(|T_\infty^1(p)|,|T_\infty^2(p)|)$
have the same joint distribution. It is not hard to see, 
that if $|T_\infty^1(p)|=k,$ then the conditional distribution of
$H(T_\infty^1(p))$ is $\cT_k.$ By construction, if $|\tilde{T}^1(p)|=k,$ then also $H(\tilde{T}^1(p)) \sim \cT_k.$
Furthermore, it is elementary to show that if $|T_\infty^1(p)|=\infty,$ then the conditional distribution of
$H(T_\infty^1(p))$ is $\cT_\infty(p).$ 
 By construction, if $|\tilde{T}^1(p)|=\infty,$
then also $H(\tilde{T}^1(p)) \sim \cT_\infty(p).$
We can therefore conclude that
$(\tilde{T}^1(p),\tilde{T}^2(p))$ and $(T_\infty^1(p),T_\infty^2(p))$ have the same joint distribution, from which the statement
follows.
\fbox{} \\

\noindent
{\bf Remark.} The crucial part in the argument was to show that  $(|\tilde{T}^1(p)|,|\tilde{T}^2(p)|)$ and $(|T_\infty^1(p)|,|T_\infty^2(p)|)$ had the same joint distribution. From this it followed quite easily that also $(\tilde{T}^1(p),\tilde{T}^2(p))$ and $(T_\infty^1(p),T_\infty^2(p))$ had the same joint distribution.
Similar situations will occur throughout the paper.

\medskip

We can now prove Theorem \ref{thmmain2} for $d=2$.

\medskip\noindent
{\bf Proof of Theorem \ref{thmmain2} for $d=2$.}
We will prove the statement through induction in $k,$ so we start by noting that the statement is trivial for $k=0,1.$
Fix $k\geq 1,$ and assume that the statement holds for any $l\leq k.$

Let $L_{k+1}$ be a random variable such that
\[
\P(L_{k+1}=l)=\left\{\begin{array}{ccc}
2 \frac{c_l c_{k-l}}{c_{k+1}} & \textrm{if} & 0\leq l<k-l, \\
\frac{c_l c_{k-l}}{c_{k+1}} & \textrm{if} & l=k-l \\
0 & \textrm{otherwise.}
\end{array} \right.
\]
When $d=2,$ the numbers $c_k$ are the Catalan numbers. It is an elementary exercise, 
to show that the above probabilities sum to one.
Let $L^*$ be a random variable such that $\P(L^*=l)=2p \eta_l(p)$ for any $0\leq l<\infty$ and $\P(L^*=\infty)=p\eta_\infty(p).$
These probabilities sum to one as explained when we defined $T^*(p).$
We observe that by (\ref{eqn5}),
\[
2p\eta_l(p)=2c_lp^{l+1}(1-p)^{l+1}
\leq 2 \frac{c_l c_{k-l}}{c_{k+1}}.
\]
Using this, it is not hard to see that we can in fact couple $L_{k+1}$ and $L^*$ such that $\P(L_{k+1} \leq L^*)=1.$

We will construct $T_{k+1}$ and $\tilde{T}(p)$ so that $T_{k+1}\sim \cT_{k+1},$ 
$\tilde{T}(p)\sim \cT_\infty (p)$ and $T_{k+1} \subset \tilde{T}(p).$ 
Informally, the tree $T_{k+1}$ is constructed by taking a root, and then 
attaching two subtrees
onto it. The size of the smallest of these subtrees is $L_{k+1},$ while the other 
will have size $k-L_{k+1}.$ By using $L^*$ (coupled with $L_{k+1}$ so that 
$L_{k+1}\leq L^ *$) to simultaneously construct $\tilde{T}(p)$ we will make
sure that $T_{k+1} \subset \tilde{T}(p).$ By the use of $L_{k+1}$ and $L^*,$
it will be straightforward to check, using Lemma \ref{lemmad=21},
 that the distributions of $T_{k+1}$ and 
$\tilde{T}(p)$ are as claimed.

In order to give the formal construction, we consider the random variables
\[
(L_{k+1},L^*,T_{0,1},T_{1,1},\ldots,T_{k,1},T_{\infty,1}(p),T_{0,2},T_{1,2},\ldots,T_{k,2},T_{\infty,2}(p),X,(T_{l,3})_{l \geq 0}),
\]
on a common probability space. The five groups \sloppy $(L_{k+1},L^*),$ $(T_{0,1},T_{1,1},\ldots,T_{k,1},T_{\infty,1}(p)),$
$(T_{0,2},T_{1,2},\ldots,T_{k,2},T_{\infty,2}(p)),$ $X$ and $(T_{l,3})_{l \geq 0}$ of random variables are independent
of each other.
Furthermore, they have the following joint distributions.
\begin{itemize}
\item $L_{k+1},L^*$ are coupled so that $L_{k+1} \leq L^*.$ 

\item For $i=1,2,$ $T_{0,i},T_{1,i},\ldots,T_{k,i},T_{\infty,i}(p)$ have marginal distributions
$T_{0,i}=\emptyset,$ $T_{l,i}\sim \cT_l$ for every $1\leq l \leq k,$ and $T_{\infty,i}(p) \sim \cT_\infty (p).$ Furthermore,
they are coupled so that
$T_{1,i} \subset \cdots \subset T_{k,i} \subset T_{\infty,i}(p).$ Such a coupling exists by Theorem \ref{thmfinite}
and the induction hypothesis.

\item $X\in\{1,2\}$ is such that $\P(X=1)=\P(X=2)=1/2.$

\item The elements of the sequence $(T_{l,3})_{l \geq 0}$ have marginal
distributions $T_{0,3}=\emptyset,$ $T_{l,3}\sim \cT_l$ for every $1\leq l <\infty.$ Furthermore,
they are coupled so that
$T_{0,3} \subset T_{1,3}\subset \cdots \subset T_{l,3}\subset \cdots.$ This is possible by Theorem \ref{thmfinite}.

\end{itemize}
On this probability space we construct $T_{k+1}$ and $\tilde{T}(p)$ as follows. Let
\[
V(T_{k+1})=\left\{\begin{array}{lcc}
\{o\}\bigcup_{u \in V(T_{L_{k+1},3})}\{(X,u)\}\bigcup_{v\in V(T_{k-L_{k+1},1})}\{(3-X,v)\} & \textrm{if} & L^*<\infty, \\
\{o\}\bigcup_{u \in V(T_{L_{k+1},1})}\{(X,u)\}\bigcup_{v\in V(T_{k-L_{k+1},2})}\{(3-X,v)\} & \textrm{if} & L^*=\infty,
\end{array} \right.
\]
and
\[
V(\tilde{T}(p))=\left\{\begin{array}{lcc}
\{o\}\bigcup_{u \in V(T_{L^*,3})}\{(X,u)\}\bigcup_{v\in V(T_{\infty,1}(p))}\{(3-X,v)\} & \textrm{if} & L^*<\infty, \\
\{o\}\bigcup_{u \in V(T_{\infty,1}(p))}\{(X,u)\}\bigcup_{v\in V(T_{\infty,2}(p))}\{(3-X,v)\} & \textrm{if} & L^*=\infty.
\end{array} \right.
\]
Informally, we use $X$ to determine which of the children of the root will be 
given the smallest number of offspring. If $L^*<\infty,$ then we attach finite 
subtrees to this child for both $T_{k+1}$ and $\tilde{T}(p),$ while if 
$L^*=\infty,$ we attach an infinite subtree to $\tilde{T}(p)$ and a finite to 
$T_{k+1}.$ 
We note that by construction $T_{k+1} \subset \tilde{T}(p).$ This can easily be checked case by case.

As mentioned above, the use of $L^*$ makes sure that $\tilde{T}(p)$ is constructed as in Lemma \ref{lemmad=21}. By that 
lemma, we conclude that $\tilde{T}(p) \sim \cT_\infty(p).$
It only remains to show that $T_{k+1}\sim \cT_{k+1}.$ 
It is easily checked that for $T\sim \cT_{k+1},$ $\min(|T^1|, |T^2|)$
has the same distribution as $L_{k+1}.$ From this we conclude that  
$(|T^1|, |T^2|)$ and $(|T_{k+1}^1|, |T_{k+1}^2|)$ have the same joint distribution.
Furthermore, for $i=1,2$ and conditional on the event $|T^i|=l,$ we get that 
$T^i \sim \cT_l.$ This follows as in the proof of Lemma
\ref{lemmad=21}, see also the remark thereafter. We conclude that indeed
 $T_{k+1}\sim \cT_{k+1}.$
\fbox{} \\

Recall the definition of $L^*=L^*(p)$ in the proof of Theorem \ref{thmmain2} above.
We will use our next lemma to prove Theorem \ref{thmmain}.
\begin{lemma} \label{lemmad=22}
Let $1/2\leq p_1<p_2.$ There exists a coupling of $L^*(p_1)$ and $L^*(p_2)$ such that $\P(L^*(p_1)\leq L^*(p_2))=1.$
\end{lemma}
{\bf Proof.}
Observe that for any $k<\infty,$ $2p\eta_k(p)=2c_k p^{k+1}(1-p)^{k+1},$
which is decreasing in $p$ when $p\geq1/2.$ Hence, for every $k<\infty$,
\begin{equation} \label{eqn23}
\sum_{j=0}^k 2p_2\eta_j(p_2)\leq \sum_{j=0}^k 2p_1\eta_l(p_1).
\end{equation}
The statement follows easily from (\ref{eqn23}) and the definition of $L^*(p).$
\fbox{} \\

We can now prove Theorem \ref{thmmain} in the case of $d=2$, by using Lemma \ref{lemmad=22} and Theorem \ref{thmmain2}. \\

\noindent
{\bf Proof of Theorem \ref{thmmain} in the case of $d=2$.}

\noindent
Let $((X_u,L_u^*(p_1),L_u^*(p_2), (T_{l,u})_{l \geq 0},T_{\infty,u}(p_1),T_{\infty,u}(p_2)))_{u \in V(\T^2)}$
be an i.i.d.\ collection, indexed by
$V(\T^2),$ and with the following distribution. The marginal distributions of the random variables
$(X_u,L_u^*(p_1),L_u^*(p_2), (T_{l,u})_{l \geq 0},T_{\infty,u}(p_1),T_{\infty,u}(p_2))$
are as indicated by
the notation, they are explained on multiple occasions above.
For fixed $u\in V(\T^2)$ the joint distribution is as follows:
\begin{itemize}
\item $X_u$ is independent of the other random variables. 

\item $L_u^*(p_1),L_u^*(p_2)$ are independent of the other random variables and coupled so that $\P(L_u^*(p_1)\leq L_u^*(p_2))=1.$
This is possible by Lemma \ref{lemmad=22}.

\item $(T_{l,u})_{l \geq 0},T_{\infty,u}(p_2)$ are independent of the other random variables and coupled so that
$T_{0,u}\subset T_{1,u}\subset\cdots\subset T_{l,u} \cdots \subset T_{\infty,u}(p_2).$ This is possible by
Theorems \ref{thmfinite} and \ref{thmmain2}.

\item $T_{\infty,u}(p_1)$ is independent of the other random variables.

\end{itemize}

We will use the random variables above to construct a sequence  $(S_n(p_1),S_n(p_2))_{n\geq 0}$ of pairs of trees
such that $S_n(p_1) \subset S_n(p_2)$ for every $n.$ We will then show that the limiting objects $S_\infty(p_i)$
is such that $S_\infty(p_i)\sim \cT_\infty(p_i)$ for $i=1,2$ and
$S_\infty(p_1) \subset S_\infty(p_2),$ thus proving the theorem.
The construction will be performed in steps, and to that end we use an ordering of
$V(\T^2).$ We simply let $o<(1)<(2)<(1,1)<(1,2)<(2,1)<\cdots,$ and proceed in the natural way.
Let $\U_0=\{o\}$ and $S_0(p_1)=S_0(p_2)=\{o\}.$ 

Before we give the formal construction of $(S_n(p_1),S_n(p_2))_{n\geq 0}$, let us explain the idea. 
Assume therefore that $n-1$ steps of the procedure has been performed. Then, 
$\U_{n-1}$ will be the set of leaves of $S_{n-1}(p_1)$ that eventually will have infinitely many descendants.
In fact, as we will see below, if we were to attach independent copies of trees with distribution $\cT_\infty(p_i)$ to 
$S_{n-1}(p_i)$ at all of the vertices of $\U_{n-1},$ we would get a tree which again would have distribution $\cT_\infty(p_i).$
We let  $u_n$ be the smallest vertex of $\U_{n-1}$ (in the ordering of $V(\T^2)$),
and then we use $X_{u_n}$ to pick one of the children of $u_n.$ If $L_{u_n}^*(p_1)<\infty$ then we use that 
$T_{L_{u_n}^*(p_1),{u_n}}\subset T_{L_{u_n}^*(p_2),{u_n}}$ and attach these trees to the vertex $(u_n,X_{u_n})$ in 
$S_n(p_1)$ and $S_n(p_2)$ respectively. For convenience, we abuse the notation somewhat and 
write $T_{L_{u_n}^*(p_2),{u_n}}$ for $T_{L_{u_n}^*(p_2),{u_n}}(p_2)$ when $L_{u_n}^*(p_2)=\infty.$
We also attach $(u_n,3-X_{u_n})$ to $S_n(p_i)$ and 
create $\U_n$ by removing $u_n$ from $\U_{n-1}$ and adding $(u_n,3-X_{u_n})$ 
(and thereby designating $(u_n,3-X_{u_n})$ to eventually have an infinite 
number of descendants). 
If instead $L_{u_n}^*(p_1)=\infty,$ then we attach the vertices $(u_n,X_{u_n})$ and $(u_n,3-X_{u_n})$
to both $S_n(p_1)$ and $S_n(p_2)$ and get $\U_n$ by removing $u_n$ from $\U_{n-1}$ and adding 
$(u_n,X_{u_n})$ and $(u_n,3-X_{u_n}).$
 The gain is that we then have $S_n(p_1)\subseteq S_n(p_2)$ and that attaching independent copies of trees 
with distribution $\cT_\infty(p_i)$ to 
$S_n(p_i)$ at all of the vertices of $\U_{n},$ we would get a tree which again would have distribution $\cT_\infty(p_i)$
(because of Lemma \ref{lemmad=21}).\\

\noindent
Formally, the construction at step $n \geq 1$ consists of the following:

\noindent
Let $u_n=\min\{u\in V(\T^2): u\in \U_{n-1}\},$ and for $i=1,2,$ set
\[
V(S_n(p_i))=\left\{\begin{array}{lcc}
V(S_{n-1}(p_i))\bigcup_{v\in V(T_{L_{u_n}^*(p_i),{u_n}})}\{({u_n},X_{u_n},v)\}\bigcup\{({u_n},3-X_{u_n})\} & \textrm{if} & L_{u_n}^*(p_1)<\infty \\
V(S_{n-1}(p_i))\bigcup \{({u_n},X_{u_n})\}\bigcup\{({u_n},3-X_{u_n})\} & \textrm{if} & L_{u_n}^*(p_1)=\infty,
\end{array} \right.
\]
and
\[
\U_n=\left\{\begin{array}{lcc}
(\U_{n-1}\setminus\{u_n\})\bigcup\{({u_n},3-X_{u_n})\} & \textrm{if} & L_{u_n}^*(p_1)<\infty \\
(\U_{n-1}\setminus\{u_n\})\bigcup\{({u_n},X_{u_n})\}\bigcup\{({u_n},3-X_{u_n})\} & \textrm{if} & L_{u_n}^*(p_1)=\infty.
\end{array} \right.
\]
It is elementary to check, using the itemized description 
above, that for every $n$ we have that $S_n(p_1) \subseteq S_n(p_2).$

Furthermore, for $i=1,2$ we define $\tilde{S}_n(p_i)$ by
\[
V(\tilde{S}_n(p_i))=V(S_n(p_i))\bigcup_{u\in \U_n} \bigcup_{v\in V(T_{\infty,u}(p_i))} \{(u,v)\}.
\]
Thus, we get $\tilde{S}_n(p_i)$ from $S_n(p_i)$ by attaching,
to every $u\in \U_n,$ an independent tree with distribution $\cT_\infty(p_i).$
We claim that $\tilde{S}_n(p_i) \sim \cT_\infty(p_i)$ for every $n$ (which is basically the reason for introducing them)
which we prove by induction. 
Consider therefore $\tilde{S}_1(p_1).$ We see that
$\P\left(\min\left(|\tilde{S}^1_1(p_1)|,|\tilde{S}^2_1(p_1)|\right)=k\right)=\P(L^*_o(p_1)=k)=2p_1\eta_k(p_1).$
If $L^*_o(p_1)<\infty$ then $\U_1=\{(3-X_o)\}$ while if $L^*_o(p_1)=\infty$ then $\U_1=\{(1),(2)\}$. Therefore, by attaching
independent trees with distribution $\cT_\infty(p_i)$ at all $u\in \U_1,$ we see that
$\tilde{S}_1(p_1)$ is constructed as $\tilde{T}(p_1)$ in the statement of Lemma \ref{lemmad=21}, and by that lemma
we have that $\tilde{S}_1(p_1) \sim \cT_\infty (p_1).$ 

Assume now that for some fixed $n,$ $\tilde{S}_n(p_1) \sim \cT_\infty(p_1).$ 
We construct $S_{n+1}(p_1)$ from $S_{n}(p_1)$ by performing the above construction at $u_{n+1}.$
Thus, when performing the constructions of $\tilde{S}_{n}(p_1)$ and $\tilde{S}_{n+1}(p_1)$
(from $S_{n}(p_1)$ and $S_{n+1}(p_1)$ respectively)
we can attach the same independent trees with distribution $\cT_\infty(p_1)$ at every $u\in \U_{n}\setminus\{u_{n+1}\}.$
When performing the rest of the construction of $\tilde{S}_{n+1}(p_1)$ at the children of $u_{n+1}$
that belongs to $\U_{n+1},$ we claim that $H(\tilde{S}^{u_{n+1}}_{n+1}(p_1)) \sim \cT_\infty (p_1)$ (here, $\tilde{S}^{u_{n+1}}_{n+1}(p_1)$ should be thought of
as $(\tilde{S}_{n+1}(p_1))^{u_{n+1}}$, that is, as a subtree of 
$\tilde{S}_{n+1}(p_1)).$
Therefore, we can in fact take $\tilde{S}_{n}(p_1)=\tilde{S}_{n+1}(p_1),$ and so 
we only need to check that $H(\tilde{S}^{u_{n+1}}_{n+1}(p_1)) \sim \cT_\infty (p_1).$
However, this follows as for $\tilde{S}_1(p_1)$ since  we have that
$H(\tilde{S}^{u_{n+1}}_{n+1}(p_1))$ is constructed as $\tilde{T}(p_1)$ in the statement of Lemma \ref{lemmad=21}, and by that lemma
we get that $H(\tilde{S}^{u_{n+1}}_{n+1}(p_1)) \sim \cT_\infty (p_1).$ 
The same argument shows that also $\tilde{S}_n(p_2) \sim \cT_\infty(p_2)$ for every $n.$

Define $S_\infty (p_i)$ by
\[
V(S_\infty (p_i))=\bigcup_{m=1}^\infty \bigcap_{n=m}^\infty V(\tilde{S}_n (p_i))=\bigcup_{n=1}^\infty V(S_n(p_i)),
\]
so that $S_\infty (p_1) \subset S_\infty (p_2).$ For any finite $A\subset \T^2$, let $\max(A)=\max\{v\in V(\T^2): v\in A\}$
where the maximum is taken with respect to the ordering of $V(\T^2),$ and let $N(A)=|\{u\in V(\T^2):u\leq \max(A)\}|.$
We get that for any $n\geq N(A),$
\[
\P(A\subset S_\infty (p_i))=\P(A\subset S_n (p_i))=\P(A\subset \tilde{S}_n (p_i)),
\]
and since $\tilde{S}_n(p_i)\sim \cT_\infty (p_i),$
the distribution of $S_\infty (p_i)$ equals $\cT_\infty (p_i)$ on any cylinder event,
and so we conclude that $S_\infty(p_i)\sim \cT_\infty (p_i).$
\fbox{} \\

\section{The case $d=3$} \label{secd=3}

As the title suggests, we will assume throughout this section that $d=3$, and also that $p \geq 1/3$.
Furthermore, we want to use similar notation as in Section \ref{secd=2}, and therefore we consider the
definitions of Section \ref{secd=2} void.
For instance, when we in this section refer to a tree with distribution $\cT_\infty(p),$ we are
implicitly assuming that $d=3.$

The approach of this section is similar to when $d=2$. Consider now a tree $ T_\infty \sim \cT_\infty(p)$ and the 
three subtrees $T^1_\infty, T^2_\infty,T^3_\infty.$ One of these will necessarily be infinite, while the 
other ones may be finite. Thus, one way of generating $ T_\infty$ should be to start with a root o, 
then to pick one of the children $(1),(2),(3)$ with equal probability and attach an independent copy from 
$\cT_\infty(p)$ to it. Then, we use random trees $T^*(p),T^{**}(p)$ 
(with a very particular joint distribution),
and attach these trees to the other children. This is made precise in Lemma \ref{lemmad=31}. 
In Lemma \ref{lemmad=32}, we then prove that for any $1/3 \leq p_1<p_2\leq 1,$ we can couple
$T^*(p_i),T^{**}(p_i)$ so that $|T^*(p_1)|\leq|T^*(p_2)| $ and $|T^{**}(p_1)|\leq|T^{**}(p_2)|.$
We can then use this together with Theorem \ref{thmfinite} to prove Theorem \ref{thmmainsc} which is the special case of 
Theorem \ref{thmmain} where $p_1=1/3.$ From this, we can then prove Theorem \ref{thmmain2}, and in turn 
Theorem \ref{thmmain}.

Our first aim of this section is to arrive at a result that is the analogue of Lemma \ref{lemmad=21}, but
for $d=3$. To that end, we will need two technical lemmas, Lemma \ref{lemmafdecrease} and
Lemma \ref{lemmad=3new1}. Observe that
$\eta_{\infty}(p)=p(1-(1-\eta_{\infty}(p))^3)=p(3\eta_{\infty}(p)-3\eta_{\infty}(p)^2+\eta_{\infty}(p)^3).$
It follows that
$
p\eta_{\infty}(p)^2-3p\eta_{\infty}(p)+3p-1=0,
$
from which we conclude that
\begin{equation} \label{eqnd=3inf2}
\eta_{\infty}(p)
=\frac{1}{2}\left(3-\sqrt{\frac{4}{p}-3}\right)
=\frac{3\sqrt{p}-\sqrt{4-3p}}{2\sqrt{p}}.
\end{equation}
Consider the function
\[
f(p):=\frac{\sqrt{3p}-1}{1-\sqrt{3p}(1-\eta_{\infty}(p))}=\frac{\sqrt{3p}-1}{1-\frac{\sqrt{3}}{2}(\sqrt{4-3p}-\sqrt{p})}.
\]
The reason for introducing $f(p)$ will become clear later, it will play a crucial part in this section. We can now state the first
of the two previously announced lemmas.
\begin{lemma} \label{lemmafdecrease}
We have that $\lim_{p \downarrow 1/3}f(p)=1,$ $f(1)=\sqrt{3}-1$ and that
 $f'(p)<0$ if $1/3< p\leq 1.$ Therefore $0\leq f(p) \leq 1$ for every $p\in[1/3,1].$
\end{lemma}
{\bf Proof.}
The statement that $f(1)=\sqrt{3}-1$ is trivial.

Using the standard expansion $\sqrt{3p}=1+(3p-1)/2+O((3p-1)^2),$ and similar expressions for $\sqrt{4-3p}$ and $\sqrt{p},$
we get that
\[
\lim_{p\downarrow 1/3}f(p)=\lim_{p\downarrow 1/3}\frac{(3p-1)/2+O((3p-1)^2)}{(3p-1)/2+O((3p-1)^2)}=1,
\]
proving the first part of the statement. Furthermore,
\[
f'(p)=\frac
{\frac{3}{2\sqrt{3p}}(1-\frac{\sqrt{3}}{2}(\sqrt{4-3p}-\sqrt{p}))
-(\sqrt{3p}-1)(-\frac{\sqrt{3}}{2}(\frac{-3}{2\sqrt{4-3p}}-\frac{1}{2\sqrt{p}}))}
{(1-\frac{\sqrt{3}}{2}(\sqrt{4-3p}-\sqrt{p}))^2},
\]
so that $f'(p)< 0$ iff
\begin{eqnarray*}
\lefteqn{0>\frac{1}{\sqrt{p}}(1-\frac{\sqrt{3}}{2}(\sqrt{4-3p}-\sqrt{p}))
-(\sqrt{3p}-1)(\frac{3}{2\sqrt{4-3p}}+\frac{1}{2\sqrt{p}})} \\
& & =\frac{2\sqrt{4-3p}-\sqrt{3}(4-3p-\sqrt{4-3p}\sqrt{p})
-(\sqrt{3p}-1)(3\sqrt{p}+\sqrt{4-3p})}{2\sqrt{4-3p}\sqrt{p}} \\
& & =\frac{3\sqrt{4-3p}-4\sqrt{3}
+3\sqrt{p}}{2\sqrt{4-3p}\sqrt{p}}.\\
\end{eqnarray*}
Therefore, we need to show that
$
3\sqrt{4-3p}<\sqrt{3}(4-\sqrt{3p}).
$
A straightforward calculation shows that this condition is the same as
$
9p^2-6p+1>0,
$
which is easily seen to be true for $p>1/3.$
\fbox{}\\

We will now give the construction of the pair of random trees $(T^*(p),T^{**}(p))$ mentioned above. 
Let $(T_{k,1})_{k\geq 0},(T_{k,2})_{k\geq 0}, T_{\infty,1}(p), T_{\infty,2}(p),U_1,U_2$ 
be independent random variables with the following marginal distributions:
\[
U_1,U_2\sim U[0,1], \ T_{k,i}\sim \cT_k \textrm{ and } T_{\infty,i}(p)\sim \cT_{\infty}(p) \textrm{ for } i=1,2.
\]
We let, for every $k<\infty,$
\[
V(T^*(p))=\left\{\begin{array}{ccl}
T_{k,1} & \textrm{if} & \sum_{l=0}^{k-1}\sqrt{3p}\eta_{l}(p)<U_1 \leq \sum_{l=0}^{k}\sqrt{3p}\eta_{l}(p),\ \  \\
T_{\infty,1}(p) & \textrm{if} & \sqrt{3p}(1-\eta_\infty(p))<U_1.
\end{array}\right.
\]
Furthermore, we define the conditional distribution of $T^{**}(p),$ given $T^*(p)$ by, for every $k<\infty,$
\[
V(T^{**}(p))=\left\{\begin{array}{cl}
T_{k,2} & \begin{array}{cl}
\textrm{if} & \sum_{l=0}^{k-1}\sqrt{3p}\eta_{l}(p)<U_2 \leq \sum_{l=0}^{k}\sqrt{3p}\eta_{l}(p),
 \textrm{ and } |T^*(p)|<\infty, \\
\textrm{or if} & \sum_{l=0}^{k-1}\sqrt{3p}\eta_{l}(p)f(p)
<U_2 \leq \sum_{l=0}^{k}\sqrt{3p}\eta_{l}(p)f(p),  \textrm{ and } |T^*(p)|=\infty
\end{array} \\
& \\
T_{\infty,2}(p) &\begin{array}{cl}
\textrm{if} & \sqrt{3p}(1-\eta_\infty(p))<U_2, \textrm{ and } |T^*(p)|<\infty,\\
\textrm{or if} & 1-\frac{p\eta_\infty(p)^2}{1-\sqrt{3p}(1-\eta_\infty(p))}
<U_2  \ \  \textrm{ and } |T^*(p)|=\infty.
\end{array} \\
\end{array}\right.
\]

We have the following lemma.
\begin{lemma}\label{lemmad=3new1}

The pair $(T^*(p),T^{**}(p))$ is well defined.

\end{lemma}
{\bf Proof.}
We need to show that all the claimed probabilities are nonnegative, and that the appropriate sums add to one.
We have that $\sum_{l=0}^\infty\P(|T^*(p)|=l)=\sum_{l=0}^{\infty}\sqrt{3p}\eta_{l}(p)=
\sum_{l=0}^{\infty}\sqrt{3}c_lp^{l+1/2}(1-p)^{2l+1},$ by (\ref{eqn10}). This sum is easily seen to be
maximized when $p=1/3,$ when it takes the value 1. Since $\sum_{l=0}^{\infty}\sqrt{3p}\eta_{l}(p)=
\sqrt{3p}(1-\eta_\infty(p)),$
it follows that $T^*(p)$ is well defined. This also proves that the conditional
distribution of $T^{**}(p),$ given the event that $|T^{*}(p)|<\infty,$ is well defined.

It remains to prove that also the conditional
distribution of $T^{**}(p),$ given the event that $|T^{*}(p)|=\infty,$ is well defined.
It follows from Lemma \ref{lemmafdecrease} and the first paragraph of this proof,
that $\sum_{l=0}^\infty \P(|T^{**}(p)|=l |\, |T^{*}(p)|=\infty)
=\sum_{l=0}^{\infty}\sqrt{3p}\eta_{l}(p)f(p)
\leq f(p)\leq 1.$ Furthermore, by the calculation leading up to (\ref{eqnd=3inf2}),
\begin{eqnarray*}
\lefteqn{\sum_{l=0}^{\infty}\sqrt{3p}\eta_l(p) f(p)
 =\sqrt{3p}(1-\eta_{\infty}(p)) \frac{\sqrt{3p}-1}{1-\sqrt{3p}(1-\eta_{\infty}(p))}} \\
& & =\frac{1-\sqrt{3p}(1-\eta_{\infty}(p))-(3p\eta_{\infty}(p)-3p+1)}{1-\sqrt{3p}(1-\eta_{\infty}(p))}
=1-\frac{p\eta_{\infty}(p)^2}{1-\sqrt{3p}(1-\eta_{\infty}(p))}.
\end{eqnarray*}
We conclude that $\sum_{l=0}^\infty \P(|T^{**}(p)|=l |\, |T^{*}(p)|=\infty)
+\P(|T^{**}(p)|=\infty |\, |T^{*}(p)|=\infty)=1$ and that all the terms of this sum are nonnegative.
\fbox{}\\

We let 
the joint distribution of the pair $(T^*(p),T^{**}(p))$ be denoted by $\cT^{*,**}(p).$
We can now present the analogue of Lemma \ref{lemmad=21}.
\begin{lemma} \label{lemmad=31}
Let $X_1,X_2,X_3,T^*(p), T^{**}(p),T_\infty(p)$ be random variables such that
\begin{itemize}
\item $X_1,X_2,X_3\in\{1,2,3\}$ are independent of the other random variables, and $(X_1,X_2,X_3)$ is a uniformly
chosen permutation of $(1,2,3).$

\item $T^*(p), T^{**}(p)$ are independent of the other random variables, and $(T^*(p),T^{**}(p))\sim\cT^{*,**}(p).$

\item $T_\infty(p)$ is independent of the other random variables, and $T_\infty(p)\sim \cT_\infty(p).$
\end{itemize}
Define the tree $\tilde{T}(p)$ by letting
\[
V(\tilde{T}(p))=\{o\}\bigcup_{u\in V(T^*(p))}\{(X_1,u)\}\bigcup_{v\in V(T^{**}(p))}\{(X_2,v)\}\bigcup_{w\in V(T_\infty(p))}\{(X_3,w)\}.
\]
We have that $\tilde{T}(p)\sim \cT_\infty(p).$
\end{lemma}
{\bf Proof.}
We start by showing that $(|\tilde{T}^1(p)|, |\tilde{T}^2(p)|, |\tilde{T}^3(p)|)$ and 
$(|T_\infty^1(p)|, |T_\infty^2(p)|, |T_\infty^3(p)|)$
have the same joint distribution.

Let $(i_1,i_2,i_3)$ be any permutation of $(1,2,3).$ We have that for any $k_1,k_2<\infty$ and $T(p)\sim \cT(p),$
\[
\begin{array}{l}
{\mathbb P}(|T^{i_1}(p)|=\infty, |T^{i_2}(p)|=k_1,|T^{i_3}(p)|=k_2 \, | \, |T(p)|=\infty )  =
   \frac{p\eta_{\infty}(p)\eta_{k_1}(p)\eta_{k_2}(p)}{\eta_{\infty}(p)}=p\eta_{k_1}(p)\eta_{k_2}(p), \\
{\mathbb P}(|T^{i_1}(p)|=\infty, |T^{i_2}(p)|=\infty,|T^{i_3}(p)|=k_1 \, | \, |T(p)|=\infty)  =
   \frac{p\eta_{\infty}(p)^2\eta_{k_1}(p)}{\eta_{\infty}(p)}=p\eta_{\infty}(p)\eta_{k_1}(p), \\
{\mathbb P}(|T^{i_1}(p)|=\infty, |T^{i_2}(p)|=\infty,|T^{i_3}(p)|=\infty \, | \, |T(p)|=\infty)  =
   \frac{p\eta_{\infty}(p)^3}{\eta_{\infty}(p)}=p\eta_{\infty}(p)^2.
\end{array}
\]
Furthermore, we have that
\begin{eqnarray*}
\lefteqn{{\mathbb P}(|\tilde{T}^{i_1}(p)|=\infty, |\tilde{T}^{i_2}(p)|=k_1, |\tilde{T}^{i_3}(p)|=k_2)
={\mathbb P}(X_1=i_2,X_2=i_3, |T^{*}(p)|=k_1,|T^{**}(p)|=k_2)} \\
& & +{\mathbb P}(X_1=i_3,X_2=i_2, |T^{*}(p)|=k_2,|T^{**}(p)|=k_1)
 =\frac{\sqrt{3p}\eta_{k_1}(p)\sqrt{3p}\eta_{k_2}(p)}{3}
=p\eta_{k_1}(p)\eta_{k_2}(p).
\end{eqnarray*}
We also see that
\begin{eqnarray*}
\lefteqn{{\mathbb P}(|\tilde{T}^{i_1}(p)|=\infty, |\tilde{T}^{i_2}(p)|=\infty, |\tilde{T}^{i_3}(p)|=\infty)} \\
& & ={\mathbb P}(|T^{*}|=\infty,|T^{**}|=\infty)
 =(1-\sqrt{3p}(1-\eta_{\infty}(p)))\frac{p\eta_{\infty}(p)^2}{1-\sqrt{3p}(1-\eta_{\infty}(p))}=p\eta_{\infty}(p)^2.
\end{eqnarray*}
Finally, 
\begin{eqnarray*}
\lefteqn{{\mathbb P}(|\tilde{T}^{i_1}(p)|=\infty, |\tilde{T}^{i_2}(p)|=\infty, |\tilde{T}^{i_3}(p)|=k_1)}\\
& & =\P(X_{1}=i_3, |T^{*}(p)|=k_1,|T^{**}(p)|=\infty)
+\P(X_{2}=i_3, |T^{*}(p)|=\infty,|T^{**}(p)|=k_1)\\
& & =\frac{1}{3}\left(\sqrt{3p}\eta_{k_1}(p)(1-\sqrt{3p}(1-\eta_\infty(p)))
+(1-\sqrt{3p}(1-\eta_\infty(p)))\sqrt{3p}\eta_{k_1}(p)f(p)\right)\\
& & =\frac{\sqrt{3p}\eta_{k_1}(p)}{3}(1-\sqrt{3p}(1-\eta_\infty(p)))(1+f(p))=p\eta_{\infty}(p)\eta_{k_1}(p),
\end{eqnarray*}
by the definition of $f(p).$
The conclusion that $(\tilde{T}^1(p), \tilde{T}^2(p), \tilde{T}^3(p))$ and $(T_\infty^1(p), T_\infty^2(p), T_\infty^3(p))$
have the same joint distribution follows as in Lemma \ref{lemmad=21}, see also the remark thereafter.
\fbox{}\\

\medskip

Define the distribution of the pair of random variables $(L^*(p),L^{**}(p))$ by letting
$\P(L^*(p)=k,L^{**}(p)=l)=\P(|T^*(p)|=k,|T^{**}(p)|=l)$ for every $0\leq l,k\leq \infty.$ Note that we allow
both $k$ and $l$ to be infinite.
Our next lemma is the analogue of Lemma \ref{lemmad=22}
for $d=3$. It is here that the function $f(p)$ and Lemma \ref{lemmafdecrease} comes to full use.
\begin{lemma} \label{lemmad=32}
Let $1/3\leq p_1<p_2 \leq 1.$ There exists a coupling of $(L^*(p_1),L^{**}(p_1))$ and
$(L^*(p_2),L^{**}(p_2))$ such that
\[
\P(L^*(p_1) \leq L^*(p_2), L^{**}(p_1)\leq L^{**}(p_2))=1.
\]
\end{lemma}
{\bf Proof.}
We will prove the statement by considering the construction of $(T^*(p),T^{**}(p))$ in Lemma \ref{lemmad=3new1} for
$p_1$ and $p_2$ simultaneously. We will show that this results in $|T^*(p_1)|\leq |T^*(p_2)|$ and
$|T^{**}(p_1)|\leq |T^{**}(p_2)|,$ and then we will simply let
$L^*(p_i)=|T^*(p_i)|$ and $L^{**}(p_i)=|T^{**}(p_i)|$ for $i=1,2.$

Observe that by (\ref{eqn10}), we have that for any $l<\infty,$ $\sqrt{3p}\eta_l(p)
=\sqrt{3}c_lp^{l+1/2}(1-p)^{2l+1}, $ which is clearly decreasing in $p$ for $p\geq 1/3.$ 
Consider the use of the random variable $U_1$ in the construction of Lemma \ref{lemmad=3new1}. 
It follows immediately, that by using the same random variable $U_1,$ for both 
$T^*(p_1)$ and $T^*(p_2),$ the construction yields $|T^*(p_1)|\leq |T^*(p_2)|.$

We now need to consider three cases depending on the values of $|T^*(p_1)|$ and $|T^*(p_2)|.$

\noindent
{\bf Case 1} ($|T^*(p_1)| \leq |T^*(p_2)|<\infty$): This case is treated exactly as when coupling $|T^*(p_1)|$ and $|T^*(p_2)|.$
We conclude that by using the same random variable $U_2$ for both $T^{**}(p_1)$ and $T^{**}(p_2),$ we can couple $|T^{**}(p_1)|$ and $|T^{**}(p_2)|$
so that $|T^{**}(p_1)|\leq |T^{**}(p_2)|.$

\noindent
{\bf Case 2} ($|T^*(p_1)| <\infty$ and $|T^*(p_2)| =\infty$): By Lemma \ref{lemmafdecrease}, $f(p_2)\leq 1,$ and so,
$\sum_{l=0}^k f(p_2)\sqrt{3p_2}\eta_l(p_2)\leq \sum_{l=0}^k \sqrt{3p_2}\eta_l(p_2)\leq \sum_{l=0}^k \sqrt{3p_1}\eta_l(p_1).$
As above, by using the same random variable $U_2$ for both $T^{**}(p_1)$ and $T^{**}(p_2),$ the construction yields $|T^{**}(p_1)|\leq |T^{**}(p_2)|.$

\noindent
{\bf Case 3} ($|T^*(p_1)|=|T^*(p_2)| =\infty$): By Lemma \ref{lemmafdecrease}, $f(p_2)\leq f(p_1),$ and so,
$\sum_{l=0}^k f(p_2)\sqrt{3p_2}\eta_l(p_2)\leq \sum_{l=0}^k f(p_1)\sqrt{3p_1}\eta_l(p_1).$
Again, by using the same random variable $U_2$ for both $T^{**}(p_1)$ and $T^{**}(p_2),$ the construction yields $|T^{**}(p_1)|\leq |T^{**}(p_2)|.$
\fbox{}\\
\medskip

We can now prove Theorem \ref{thmmain} in the case of $d=3$, by using Lemma \ref{lemmad=32} and Theorem \ref{thmfinite}.
In the case $d=2,$ we could use Theorem \ref{thmmain2} to prove Theorem \ref{thmmain}. Here, we have to do things
slightly differently, since we do not (yet) have a version of Theorem \ref{thmmain2}  for $d=3.$
As mentioned, we will therefore start by proving the special case $p_1=1/3.$ 
\begin{theorem} \label{thmmainsc}
Let $d=3$ and $1/3<p\leq 1.$ There exists a coupling of $T_\infty(1/3)$ and $T_\infty(p)$
(where $T_\infty(1/3) \sim \cT_\infty(1/3)$ and  $T_\infty(p) \sim \cT_\infty(p)$) such that 
\[
\P(T_\infty(1/3) \subset T_\infty(p))=1.
\]
\end{theorem}

\noindent
{\bf Proof.}
It will be convenient in what follows to set $p_1=1/3$ and $p_2=p.$ We will however only use the notation $p_1$
implicitly in formulas involving $p_i$ for $i=1,2$ and write $1/3$ in all other places.

Let
\begin{align} \label{eqn8}
((X_{1,u},X_{2,u},X_{3,u}),(L_u^*(1/3),L_u^{**}(1/3),L_u^*(p_2),  L_u^{**}(p_2)),\\
(T_{l,1,u})_{l\geq 0},(T_{l,2,u})_{l\geq 0},T_{\infty,u}(1/3),  T_{\infty,u}(p_2) )_{u\in V(\T^3)} \nonumber
\end{align}
be an i.i.d.\ collection, indexed by $V(\T^3).$  Furthermore, for any $u\in V(\T^3),$
\begin{enumerate}[(i)]
\item The random variables $(X_{1,u},X_{2,u},X_{3,u})$ have the joint distribution as described in the statement of Lemma \ref{lemmad=31}, and they are independent of all the other random variables.

\item $(L_u^*(1/3),L_u^{**}(1/3),L_u^*(p_2),L_u^{**}(p_2))$ are independent of the other random variables, the 
pairs $(L_u^*(1/3),L_u^{**}(1/3))$ and $(L_u^*(p_2),L_u^{**}(p_2))$ have joint distributions as in their constructions,
and they are coupled so that
$L_u^*(1/3) \leq L_u^*(p_2)$ and $L_u^{**}(1/3)\leq L_u^{**}(p_2).$ This is possible by Lemma \ref{lemmad=32}.

\item For any $u,$ $T_{l,1,u}\sim \cT_l.$ Furthermore, the collection
$(T_{l,1,u})_{l\geq 0},$ is independent of the other random variables and
coupled so that $T_{0,1,u}\subset T_{1,1,u}\subset \cdots.$
This is possible by Theorem \ref{thmfinite}.

\item For any $u,$ $T_{l,2,u}\sim \cT_l.$ Furthermore, the collection
$(T_{l,2,u})_{l\geq 0},$ is independent of the other random variables and
coupled so that $T_{0,2,u}\subset T_{1,2,u}\subset \cdots.$ Again, this uses Theorem \ref{thmfinite}.

\item $T_{\infty,u}(1/3)\sim \cT_{\infty,u}(1/3)$ and $T_{\infty,u}(p_2)\sim \cT_{\infty,u}(p_2)$ 
are independent of each other and all the other random variables.

\end{enumerate}

We take the ordering of $V(\T^3)$ to be the natural one, i.e. we let $o<(1)<(2)<(3)<(1,1)<(1,2)\cdots.$
Let $V_n(\T^3)$ be the set that consist of the $n$ first elements in the ordering of $V(\T^3).$

Before we give the formal construction, let us briefly explain the idea. Similar to when $d=2,$ we will use the random 
variables of (\ref{eqn8}), to construct a sequence $(S_n(1/3),S_n(p_2))_{n \geq 0}$ of pairs of trees such that 
$S_n(1/3) \subseteq S_n(p_2)$ for every $n$. Of course, here we have one more child to deal with.
The main difference is that when $d=2,$
we could divide the construction into cases depending on $L^*(p_1)$ since we had Theorem \ref{thmmain2}
at our disposal. In doing this, we made sure that both $S_n(p_1)$ and $S_n(p_2)$ were constructed in the, for us, 
appropriate way. Here, we have to divide the analogous construction into cases depending on $L^*(p_2)$
and $L^{**}(p_2)$ (i.e. we use $p_2$ instead of $p_1$) and in doing that,
we can show that the limit $S_\infty(p_2)$ has distribution $\cT_\infty(p_2).$
However, in the absence of a version of Theorem \ref{thmmain2} for $d=3,$ we will have to work a bit harder when 
it comes to $S_n(1/3).$ In fact, we will {\em not} have that $S_\infty(1/3) \sim \cT_\infty(1/3).$
Instead, for every $n,$ we will use $S_n(1/3)$ to construct yet another random tree  
$\bar{S}_n(1/3)\sim \cT_\infty(1/3)$ such that for every $n,$
$\bar{S}_n(1/3)\bigcap V_n(\T^3) \subset S_\infty(p_2).$ The statement will then follow.

 Let $\U_0=\{o\},$ and $S_0(1/3)=S_0(p_2)=\{o\}.$ We assume that $S_{n-1}(1/3),$
$S_{n-1}(p_2)$ and $\U_{n-1}$ has been constructed. Step $n\geq 1$ consists of the following:
\noindent
Let $u_n=\min\{u\in V(\T^3):u\in \U_{n-1}\},$ and for $i=1,2$ let $V(S_n(p_i))$ be equal to
\[
\begin{array}{ll}
& \textrm{ if } \\
  \begin{array}{l}
    V(S_{n-1}(p_i))\bigcup_{v\in V(T_{L^*_{u_n}(p_i),1,u_n})}\{(u_n,X_{1,u_n},v)\} \\
    \ \ \bigcup_{w\in V(T_{L^{**}_{u_n}(p_i),2,u_n})}\{(u_n,X_{2,u_n},w)\}\bigcup\{(u_n,X_{3,u_n})\}
  \end{array}
&  L^*_{u_n}(p_2),L^{**}_{u_n}(p_2)<\infty\\
& \\
  \begin{array}{l}
V(S_{n-1}(p_i))\bigcup_{v\in V(T_{L^*_{u_n}(p_i),1,u_n})}\{(u_n,X_{1,u_n},v)\} \\
\ \ \bigcup \{(u_n,X_{2,u_n})\}\bigcup\{(u_n,X_{3,u_n})\}
  \end{array}
& L^*_{u_n}(p_2)<\infty,L^{**}_{u_n}(p_2)=\infty\\
& \\
  \begin{array}{l}
V(S_{n-1}(p_i))\bigcup\{(u_n,X_{1,u_n})\} \\
\ \ \bigcup_{w\in V(T_{L^{**}_{u_n}(p_i),2,u_n})}\{(u_n,X_{2,u_n},w)\}\bigcup\{(u_n,X_{3,u_n})\}
  \end{array}
&  L^*_{u_n}(p_2)=\infty,L^{**}_{u_n}(p_2)<\infty\\
& \\
V(S_{n-1}(p_i))\bigcup\{(u_n,X_{1,u_n})\} \bigcup\{(u_n,X_{2,u_n})\}\bigcup\{(u_n,X_{3,u_n})\} &  L^*_{u_n}(p_2)=L^{**}_{u_n}(p_2)=\infty,\\
\end{array}
\]
and let
\[
\U_n=\left\{\begin{array}{lcl}
\U_{n-1}\setminus\{u_n\}\bigcup \{(u_n,X_{3,u_n})\}  & \textrm{if} & L^*_{u_n}(p_2),L^{**}_{u_n}(p_2)<\infty\\
\U_{n-1}\setminus\{u_n\}\bigcup \{(u_n,X_{2,u_n})\}\bigcup\{(u_n,X_{3,u_n})\}  & \textrm{if} &L^*_{u_n}(p_2)<\infty,L^{**}_{u_n}(p_2)=\infty \\
\U_{n-1}\setminus\{u_n\}\bigcup \{(u_n,X_{1,u_n})\}\bigcup\{(u_n,X_{3,u_n})\}  & \textrm{if} &L^*_{u_n}(p_2)=\infty,L^{**}_{u_n}(p_2)<\infty \\
\U_{n-1}\setminus\{u_n\}\bigcup \{(u_n,X_{1,u_n})\}\bigcup\{(u_n,X_{2,u_n})\}\bigcup\{(u_n,X_{3,u_n})\}
                  & \textrm{if} &L^*_{u_n}(p_2)=L^{**}_{u_n}(p_2)=\infty. \\
\end{array}
\right.
\]
Here, we abuse notation in that 
$T_{L^*_{u_n}(p_2),1,u_n}=T_{L^*_{u_n}(p_2),1,u_n}(p_2)$ 
whenever $L^*_{u_n}(p_2)=\infty$
and similarly for $L^{**}_{u_n}(p_2).$
As mentioned above, the conditions are in terms of $L^*_{u_n}(p_2)$ and $L_{u_n}^{**}(p_2),$ while in $d=2,$ the corresponding
conditions were in terms of $L^*_{u_n}(p_1).$ If we would have had a version of Theorem \ref{thmmain2} for $d=3,$
we could have coupled the sequence $(T_{l,1,u})_{l\geq 0}$ with another random tree
$T_{\infty,1,u}(p_2) \sim \cT_{\infty}(p_2)$ such that $T_{l,1,u}\subset \cdots \subset T_{\infty,1,u}(p_2).$
Then, much as when $d=2,$ we could have divided the construction into cases depending on 
$L^*_{u_n}(p_1)$ and $L_{u_n}^{**}(p_1),$ and proceeded analogously. The effect of this change in approach is described 
and dealt with below. 
Note also that by construction, $S_n(1/3)\subseteq S_n(p_2)$
for every $n,$ this can easily be checked case by case using (ii) above.

For $i=1,2$ we define $\tilde{S}_n(p_i)$ by
\begin{equation} \label{eqn2}
V(\tilde{S}_n(p_i))=V(S_n(p_i))\bigcup_{u\in\U_n}\bigcup_{v\in V(T_{\infty,u}(p_i))}\{(u,v)\}.
\end{equation}
As when $d=2,$ we want to show that $\tilde{S}_n(p_2) \sim \cT_\infty(p_2)$ for every $n.$
Consider first $\tilde{S}_1(p_2)$. By the use of the random variables $L^*_{o}(p_2),L^{**}_{o}(p_2),$
we see that $\tilde{S}_1(p_2)$ is constructed as $\tilde{T}(p_2)$ in Lemma \ref{lemmad=31}. Therefore,
$\tilde{S}_1(p_2) \sim \cT_\infty(p_2).$ 

Assume now that for some fixed $n,$ $\tilde{S}_n(p_2) \sim \cT_\infty(p_2).$ 
We construct $S_{n+1}(p_2)$ from $S_{n}(p_2)$ by performing the above construction at $u_{n+1}.$
Thus, when performing the constructions of $\tilde{S}_{n}(p_2)$ and $\tilde{S}_{n+1}(p_2)$
(from $S_{n}(p_2)$ and $S_{n+1}(p_2)$ respectively)
we can attach the same independent trees with distribution $\cT_\infty(p_2)$ at every $u\in \U_{n}\setminus\{u_{n+1}\}.$
When performing the rest of the construction of $\tilde{S}_{n+1}(p_1)$ at the children of $u_{n+1}$
that belongs to $\U_{n+1},$ we claim that $H(\tilde{S}^{u_{n+1}}_{n+1}(p_1)) \sim \cT_\infty (p_1).$ 
Therefore, we can in fact take $\tilde{S}_{n}(p_1)=\tilde{S}_{n+1}(p_1),$ and so 
we only need to check that $H(\tilde{S}^{u_{n+1}}_{n+1}(p_1)) \sim \cT_\infty (p_1).$
However, this follows as above since  we have that
$H(\tilde{S}^{u_{n+1}}_{n+1}(p_2))$ is constructed as $\tilde{T}(p_2)$ in Lemma \ref{lemmad=31}, and by that lemma
$H(\tilde{S}^{u_{n+1}}_{n+1}(p_2)) \sim \cT_\infty (p_2).$

By defining $S_\infty(p_2)$ through
\[
V(S_\infty(p_2))=\bigcup_{m=0}^\infty \bigcap_{n=m}^\infty V(\tilde{S}_n(p_2))=\bigcup_{n=0}^\infty V(S_n(p_2)),
\]
we get that $S_\infty(p_2)\sim \cT_\infty(p_2)$ exactly as when $d=2.$

However, the tree $\tilde{S}_n(1/3)$ is {\em not} distributed in accordance with Lemma \ref{lemmad=31}.
It is in fact "too big" and therefore $\tilde{S}_n(1/3)$ does not have distribution $\cT_\infty(1/3).$ 
To see this, consider $\tilde{S}_1(1/3)$ and assume that
$L_o^*(1/3)<\infty,L_o^*(p_2)=\infty$ while 
$L_o^{**}(1/3)<\infty,L_o^{**}(p_2)<\infty$
so that $\U_1=\{(X_{1,o}),(X_{3,o})\}.$
Since $(X_{1,o})\in \U_1,$ we have by construction that 
$H(\tilde{S}_1^{(X_{1,o})}(1/3))=T_{\infty,(X_{1,o})}(1/3)\sim \cT_\infty(1/3).$
However, in order for $\tilde{S}_1(1/3)$ to be constructed as in Lemma \ref{lemmad=31},
we should have let $H(\tilde{S}_1^{(X_{1,o})}(1/3)) \sim \cT_{L^*_o (1/3)}.$ 
This is an effect of using $L^*_{u_n}(p_2), L^{**}_{u_n}(p_2)$ in the construction (necessitated by
the absence of Theorem \ref{thmmain2}).

The strategy is to 'prune' the tree $\tilde{S}_n(1/3)$, without losing the inclusion property that we desire. 
Informally, we want to replace the trees that are too big by other trees of the correct size.
To that end,
for any $k\geq 2$ and $u=(u_1,\ldots,u_k),$ let $u^-=(u_1,\ldots,u_{k-1}).$ If $k=1$ we let $u^-=o.$
We let
\begin{align*}
\V_n = & \{v\in \bigcup_{k=0}^n \U_k: v=(v^-,X_{1,v^-}), \ L^*_{v^-}(1/3)<\infty, L^*_{v^-}(p_2)=\infty\} \\
& \ \ \bigcup \{v\in \bigcup_{k=0}^n \U_k: v=(v^-,X_{2,v^-}), \ L^{**}_{v^-}(1/3)<\infty, L^{**}_{v^-}(p_2)=\infty\}.
\end{align*}
It is convenient to think of $\V_n$ as the set of vertices that needs to be pruned. Note that in the example
of $\tilde{S}_1(1/3)$ above, $\V_1=\{(X_{1,o})\}.$
For $v\in \V_n,$ either $v=(v^-,X_{1,v^-})$ or $v=(v^-,X_{2,v^-})$ and we let $L_v=L^*_{v^-}(1/3)$
in the first case and $L_v=L^{**}_{v^-}(1/3)$ in the second. Thus, $L_v$ is the size that 
the subtree of $v$ should have been given if we had followed the construction of Lemma \ref{lemmad=31}.

We will perform the pruning in steps. Therefore, let $k=|\V_n|$ and $v_1<v_2<\cdots<v_{k}$ be the elements of $\V_n.$
Define the sequence $(\bar{S}_{n,i}(1/3))_{i= 1}^{k}$ of pruned trees in the following way.
The first subtree to be pruned is the one corresponding to $v_{k}$ i.e. $\tilde{S}^{v_{k}}_n(1/3).$
We have that $H(\tilde{S}^{v_k}_n(1/3)) \sim \cT_{\infty}(1/3)$ by the construction.
This follows
as when showing that $\tilde{S}_n(p_2) \sim \cT_{\infty}(p_2),$ and uses that no descendants of $v$ 
belongs to $\V_n.$ 

We will remove $\tilde{S}^{v_k}_n(1/3),$ and replace it by a tree of
size $L_{v_k}.$ Therefore, we extend our probability space by adding a
random tree $T_{n,v_k,L_{v_k}}$ with distribution $\cT_{L_{v_k}}$ and coupled with $\tilde{S}^{v_k}_n(1/3)$ so that
$T_{n,v_k,L_v}\subset H(\tilde{S}^{v_k}_n(1/3)).$ This is possible due to Theorem \ref{thmfinite}.
Furthermore, we can take $T_{n,{v_k},L_{v_k}}$ to be independent of every other random variable of (\ref{eqn8}) (except $L_{v_k}$),
which is associated to a vertex $w\in V(\T^3)$ for which there does not exist any  $u\in V(\T^3)$ such that
$w=({v_k},u).$
In other words, $T_{n,{v_k},L_{v_k}}$ only depends on $L_{v_k}$ and the random variables used to construct $\tilde{S}^{v_k}_n(1/3).$
The first pruning step is then
\[
V(\bar{S}_{n,1}(1/3))=\left(V(\tilde{S}(1/3))\setminus \bigcup_{u\in V(\T^3)}\{(v_k,u)\}\right)
\bigcup_{u\in T_{n,v_k,L_{v_k}}}\{(v_k,u)\}.
\]
In words, we first delete $v_k$ and all its descendants and then add the appropriate smaller tree.

We now proceed in the obvious manner, and assume therefore that we have 
performed $i$ pruning steps. 
We add to our probability space a tree 
$T_{n,v_{k-i},L_{v_{k-i}}}\sim \cT_{L_{v_{k-i}}},$ 
such that 
$T_{n,v_{k-i},L_{v_{k-i}}} \subset H(\tilde{S}^{v_{k-i}}_{n,i}(1/3)),$
which only depends on $L_{v_{k-i}}$ and the random variables used to construct
$\tilde{S}^{v_{k-i}}_{n,i}(1/3).$ Here, $\tilde{S}^{v_{k-i}}_{n,i}(1/3)$ 
should be thought of as $(\tilde{S}_{n,i}(1/3))^{v_{k-i}}$, that is, 
as a subtree of $\tilde{S}_{n,i}(1/3).$ We use similar notation below.
Set
\[
V(\bar{S}_{n,i+1}(1/3))=\left(V(\bar{S}_{n,i}(1/3))\setminus \bigcup_{u\in V(\T^3)}\{(v_{k-i},u)\}\right)
\bigcup_{u\in T_{n,v_{k-i},L_{v_{k-i}}}}\{(v_{k-i},u)\},
\]
and define $\bar{S}_{n}(1/3)$ through $V(\bar{S}_{n}(1/3))=V(\bar{S}_{n,k}(1/3)).$ 

By our construction, $\bar{S}_{n}(1/3) \sim \cT_\infty(1/3)$ for every $n.$ To see this, consider first 
$\bar{S}_{1}(1/3).$ By the construction of $\tilde{S}_1(1/3)$ and our pruning procedure, the size of the 
subtrees $\bar{S}^{X_{1,o}}_{1}(1/3),\bar{S}^{X_{2,o}}_{1}(1/3),\bar{S}^{X_{3,o}}_{1}(1/3)$
are $L_o^*(1/3),L_o^{**}(1/3)$ and $\infty$ respectively. Thus, by Lemma \ref{lemmad=31}, 
$H(\bar{S}_{1}(1/3)) \sim \cT_\infty(1/3).$ Assume that for fixed $n,$ $\bar{S}_{n}(1/3) \sim \cT_\infty(1/3).$
Consider $\tilde{S}_{n+1}(1/3),$ and assume first that $\V_{n+1}$ does not include any 
children of $u_{n+1}.$ This means that the size of the subtrees 
$\bar{S}^{X_{1,u_{n+1}}}_{n+1}(1/3),\bar{S}^{X_{2,u_{n+1}}}_{n+1}(1/3),\bar{S}^{X_{3,u_{n+1}}}_{n+1}(1/3)$ are 
$L_{u_{n+1}}^*(1/3),L_{u_{n+1}}^{**}(1/3)$ and $\infty$ so that by Lemma \ref{lemmad=31},
$H(\bar{S}^{u_{n+1}}_{n+1}(1/3)) \sim \cT_\infty(1/3).$ In case $\V_{n+1}$ does include a child of 
$u_{n+1},$ then by the first one or two steps of the pruning procedure (depending on whether there
are one or two children of $u_{n+1}$ in $\V_{n+1}$), 
$H(\bar{S}^{u_{n+1}}_{n+1}(1/3))$ has been replaced by a 
subtree which has distribution $\cT_\infty(1/3).$ The fact that the children of $u_{n+1}$ that belongs to 
$\V_{n+1}$ are the first to be addressed in the pruning procedure follows by the definition of 
$\V_{n+1}$ and the ordering of $V(\T^3).$
By continuing the pruning procedure 
simultaneously for both $\tilde{S}_{n}(1/3)$ and $\tilde{S}_{n+1}(1/3),$ we see that we can in fact 
take $\bar{S}_{n}(1/3)=\bar{S}_{n+1}(1/3).$

By the above construction and pruning procedure, we get that
\begin{equation} \label{eqn7}
\bar{S}_{n}(1/3)\bigcap V_n(\T^3) \subset S_n(1/3)\bigcap V_n(\T^3)
\subset S_n(p_2)\bigcap V_n(\T^3)=S_\infty(p_2)\bigcap V_n(\T^3).
\end{equation}

To conclude the theorem, let $\gamma_n$ be the measure on $\{0,1\}^{\T^3} \times \{0,1\}^{\T^3}$ with marginal distributions
$\cT_\infty(1/3)$ and $\cT_\infty(p_2)$ such that $\gamma_n(\xi(V_n(\T^3))\leq \eta(V_n(\T^3)))=1.$ The existence 
of $\gamma_n$ follows from \eqref{eqn7}. Here, we identify a tree $T$ and an element $\xi_T \in \{0,1\}^{\T^3}$ by letting
$\xi_T(v)=1$ iff $v\in T.$ Since $\{0,1\}^{\T^3} \times \{0,1\}^{\T^3}$ is compact, there exists a subsequential 
limiting measure $\gamma$ with marginal distributions $\cT_\infty(1/3)$ and $\cT_\infty(p_2)$ such that 
$\gamma(\xi(V(\T^3))\leq \eta(V(\T^3)))=\lim_n \gamma(\xi(V_n(\T^3))\leq \eta(V_n(\T^3)))=1.$
By Strassen's theorem, it follows that 
there exists random trees $S_{\infty}(1/3) \sim \cT_\infty(1/3),$ $S_{\infty}(p_2) \sim \cT_\infty(p_2),$ 
such that $\P(S_{\infty}(1/3)\subset S_{\infty}(p_2))=1.$
\fbox{}\\

\noindent
{\bf Proof of Theorem \ref{thmmain2} for $d=3$.} It follows from applying Theorems \ref{thmfinite} and
\ref{thmmainsc}.
\fbox{}\\

\noindent
{\bf Proof of Theorem \ref{thmmain} when $d=3.$}
The argument for $1/3< p_1 <p_2\leq 1$ is very similar to the proof of Theorem \ref{thmmainsc} and we will therefore only
address the necessary adjustments.

\begin{enumerate}
\item We change (iii) to state that:

(iii') $(T_{l,1,u})_{l\geq 0},$ and $T_{\infty,1,u}(p_2)$
are independent of the other random variables and coupled so that 
$T_{0,1,u}\subset T_{1,1,u}\subset \cdots \subset T_{\infty,1,u}(p_2).$ 
This is possible by using Theorem \ref{thmmain2}.

Here, $T_{\infty,1,u}(p_2)$ is added to (\ref{eqn8}). The coupling
exists, since for $(T_l)_{l\geq 0},$
$T_\infty(1/3)$ and $T_\infty(p_2)$ (with obvious distributions) we can couple these 
so that $T_0 \subset T_1 \subset \cdots \subset T_\infty(1/3) \subset T_\infty(p_2),$
using Theorems \ref{thmfinite} and \ref{thmmainsc}. We change (iv) similarly.

\item When constructing $S_n(p_i)$ and $\U_n,$ we change all conditions concerning $L^*(p_2)$ and $L^{**}(p_2)$
to the corresponding conditions for $L^*(p_1)$ and $L^{**}(p_1).$

\item We skip the entire pruning procedure and instead proceed as in the case $d=2.$

\end{enumerate}
\fbox{} \\

\section{Open problems} \label{secopenprob}

We present some open problems.

\begin{problem}
Is it possible to generalize the results of this paper to all $d\geq 4$?
\end{problem}
{\bf Remark:} Central to the case $d=3$ was to find the ``right'' function 
$f(p)$ that allowed us to construct the relevant couplings, i.e. 
a construction as in Lemma \ref{lemmad=31}, yielding the analogue of
Lemma \ref{lemmad=32}. Presumably, the approach of this paper could then 
work to solve the problem. However, 
already in $d=4,$ the analogue of this procedure becomes
much more complicated.

\begin{problem}
For which classes of parametrized offspring distributions can one obtain results such as in
this paper?
\end{problem}
{\bf Remark:} From \cite{LPS} we know that it is possible in the case of Poisson offspring distributions.

\begin{problem} \label{q3}
Consider bond percolation on ${\mathbb Z}^d$ with $p>p_c$. Consider the open cluster of the origin, conditioned
on being infinite, and denote a sample of such a cluster by ${\cal C}_{\infty}(p)$. Is it the case
that for any $p_c<p_1<p_2\leq 1 $ there exists a coupling of ${\cal C}_{\infty}(p_1)$ and  ${\cal C}_{\infty}(p_2)$
such that $\P( {\cal C}_{\infty}(p_1)\subset{\cal C}_{\infty}(p_2))=1$?
\end{problem}

{\bf Acknowledgements.} I would like to thank Russ Lyons for suggesting the 
problem and for
many useful discussions. I would also like to thank Indiana University for
supporting my visit there, a visit during which this project started. 
Finally, I would like to thank the anonymous referee(s) for many 
useful comments and suggestions.

\end{document}